\documentclass[a4paper,12pt]{article}
\usepackage{amscd}
\usepackage{amsmath,amsfonts,amssymb,amscd}
\usepackage{indentfirst,graphicx,epsfig}
\usepackage{graphicx,psfrag}
\input{epsf}

\newtheorem{thm}{Theorem}

\newtheorem{lem}{Lemma}

\def\qed{\hfill \nopagebreak\rule{5pt}{8pt}}
\def\pf{\noindent {\it Proof.} }

\title{\bf Note on the energy of regular graphs \footnote{Supported by NSFC No.10831001, PCSIRT and the
``973" program. } }
\author{
\small Xueliang Li, Yiyang Li and Yongtang Shi\\
\small Center for Combinatorics and LPMC-TJKLC \\
\small Nankai University, Tianjin 300071, China \\
\small lxl@nankai.edu.cn; liycldk@mail.nankai.edu.cn;
shi@nankai.edu.cn
\date{}}
\begin{document}
\maketitle
\begin{abstract}
For a simple graph $G$, the energy $\mathcal{E}(G)$ is defined as
the sum of the absolute values of all the eigenvalues of its
adjacency matrix $A(G)$. Let $n, m$, respectively, be the number of
vertices and edges of $G$. One well-known inequality is that
$\mathcal{E}(G)\leq \lambda_1+\sqrt{(n-1)(2m-\lambda_1)}$, where
$\lambda_1$ is the spectral radius. If $G$ is $k$-regular, we have
$\mathcal{E}(G)\leq k+\sqrt{k(n-1)(n-k)}$. Denote
$\mathcal{E}_0=k+\sqrt{k(n-1)(n-k)}$. Balakrishnan [{\it Linear
Algebra Appl.\/} {\bf 387} (2004) 287--295] proved that for each
$\epsilon>0$, there exist infinitely many $n$ for each of which
there exists a $k$-regular graph $G$ of order $n$ with $k< n-1$ and
$\frac{\mathcal{E}(G)}{\mathcal{E}_0}<\epsilon$, and proposed an
open problem that, given a positive integer $n\geq 3$, and
$\epsilon>0$, does there exist a $k$-regular graph $G$ of order $n$
such that $\frac{\mathcal{E}(G)}{\mathcal{E}_0}>1-\epsilon$. In this
paper, we show that for each $\epsilon>0$, there exist infinitely
many such $n$ that
$\frac{\mathcal{E}(G)}{\mathcal{E}_0}>1-\epsilon$. Moreover, we
construct another class of simpler graphs which also supports the
first assertion that
$\frac{\mathcal{E}(G)}{\mathcal{E}_0}<\epsilon$.\\[3mm]

\end{abstract}

\section{Introduction}

Let $G$ be a simple graph with $n$ vertices and $m$ edges. Denote by
$\lambda_1\geq \lambda_2\geq \cdots\geq \lambda_n$ the eigenvalues
of $G$. Note that $\lambda_1$ is called the spectral radius. The
energy of G is defined as $\mathcal{E}(G)=\sum_{i=1}^n|\lambda_i|$.
For more information on graph energy we refer to \cite{Gu, glz}, and
for terminology and notations not defined here, we refer to Bondy
and Murty \cite{BM1}.

One well-known inequality for the energy of a graph $G$ is that
$\mathcal{E}(G)\leq \lambda_1+\sqrt{(n-1)(2m-\lambda_1)}$. If $G$ is
$k$-regular, we have $\mathcal{E}(G)\leq k+\sqrt{k(n-1)(n-k)}$.
Denote $\mathcal{E}_0=k+\sqrt{k(n-1)(n-k)}$. In \cite{Ba},
Balakrishnan investigated the energy of regular graphs and proved
that for each $\epsilon>0$, there exist infinitely many $n$ for each
of which there exists a $k$-regular graph $G$ of order $n$ with
$k\leq n-1$ and $\frac{\mathcal{E}(G)}{\mathcal{E}_0}<\epsilon$. In
this paper, we construct another class of simpler graphs which also
support the above assertion. Furthermore, we show that for each
$\epsilon>0$, there exist infinitely many $n$ satisfying that there
exists a $k$-regular graph $G$ of order $n$ with $k< n-1$ and
$\frac{\mathcal{E}(G)}{\mathcal{E}_0}>1-\epsilon$, which answers the
following open problem proposed by Balakrishnan in \cite{Ba}:

{\bf Open problem.} Given a positive integer $n\geq 3$ and
$\epsilon>0$, does there exist a $k$-regular graph $G$ of order $n$
such that $\frac{\mathcal{E}(G)}{\mathcal{E}_0}>1-\epsilon$ for some
$k< n-1$?

\section{Main results}

Throughout this paper, we denote $V(G)$ the vertex set of $G$ and
$E(G)$ the edge set of $G$. Firstly, we will introduce the following
useful result given by So et al. \cite{ky}.

\begin{lem}\label{ky}
Let $G-e$ be the subgraph obtained by deleting an edge $e$ of
$E(G)$. Then $$\mathcal{E}(G)\leq \mathcal{E}(G-e)+2.$$
\end{lem}

We then formulate the following theorem by employing the above
lemma.

\begin{thm}[\cite{Ba}]
For any $\varepsilon>0$, there exist infinitely many $n$ for each of
which there exists a $k$-regular graph $G$ of order $n$ with $k<n-1$
and $\mathcal {E}(G)/\mathcal{E}_0<\varepsilon$.
\end{thm}

\pf Let $q>2$ be a positive integer. We take $q$ copies of the
complete graph $K_q$. Denote by $v_1,\ldots,v_q$ the vertices of
$K_q$ and the corresponding vertices in each copy by
$v_1[i],\ldots,v_q[i]$, for $1\leq i\leq q$. Let $G_{q^2}$ be  a
graph consisting of $q$ copies of $K_q$ and $q^2$ edges by joining
vertices $v_j[i]$ and $v_j[i+1]$, $(1\leq i<q)$, $v_j[q]$ and
$v_j[1]$ where $1\leq j\leq q$. Obviously, the graph $G_{q^2}$ is
$q+1$ regular. Employing Lemma \ref{ky}, deleting all the $q^2$
edges joining two copies of $K_q$, we have $\mathcal{E}(G_{q^2})\leq
\mathcal{E}(qK_q)+2q^2$. Thus, $\mathcal{E}(G_{q^2})\leq
2q(q-1)+2q^2$. Then, it follows that
\begin{align*}
\frac{\mathcal{E}(G_{q^2})}{\mathcal{E}_0}&\leq
\frac{4q^2-2q}{q+1+\sqrt{(q+1)(q^2-1)(q^2-q-1)}}\\
&\leq \frac{4q^2-2q}{(q^2-q-1)\sqrt{q+1}}\rightarrow 0 \mbox{ as }
q\rightarrow\infty.
\end{align*}
Thus, for any $\varepsilon>0$, when $q$ is large enough, the graph
$G_{q^2}$ satisfies the required condition. The proof is thus
complete.\qed

\begin{thm}
For any $\varepsilon>0$,  there exist infinitely many $n$ satisfying
that there exists a $k$-regular graph of order $n$ with $k<n-1$ and
$\mathcal{E}(G)/\mathcal{E}_0>1-\varepsilon$.
\end{thm}
\pf It suffices to verify an infinite sequence of graphs satisfying
the condition. To this end, we focus on the Paley graph (for details
see \cite{GR}). Let $p\geq 11$ be a prime and $p\equiv1 (mod\ 4)$.
The Paley graph $G_p$ of order $p$ has the elements of the finite
field $GF(q)$ as vertex set and two vertices are adjacent if and
only if their difference is a nonzero square in $GF(q)$. It is well
known that the Paley graph $G_p$ is a $(p-1)/2$-regular graph. And
the eigenvalues are $\frac{p-1}{2}$ (with multiplicity 1) and
$\frac{-1\pm\sqrt{p}}{2}$ (both with multiplicity $\frac{p-1}{2}$).
Consequently, we have
$$\mathcal{E}(G_p)=
\frac{p-1}{2}+\frac{-1+\sqrt{p}}{2}\cdot\frac{p-1}{2}+\frac{1+\sqrt{p}}{2}\cdot\frac{p-1}{2}
=(p-1)\frac{1+\sqrt{p}}{2}
>\frac{p^{3/2}}{2}.$$
Moreover,
$\mathcal{E}_0=\frac{p-1}{2}+\sqrt{\frac{p-1}{2}(p-1)({p-\frac{p-1}{2}})}$,
we can deduce that
$$\mathcal{E}(G_p)/\mathcal{E}_0>\frac{\frac{p^{3/2}}{2}}{\frac{p-1}{2}(\sqrt{p+1}+1)}
>\frac{\frac{p^{3/2}}{2}}{\frac{p}{2}(\sqrt{p}+2)}\rightarrow 1 \mbox{ as } p\rightarrow\infty.$$

Therefore, for any $\varepsilon>0$ and some integer $N$, if $p>N$,
it follows that $\mathcal{E}(G_p)/\mathcal{E}_0>1-\varepsilon$. The
theorem is thus proved. \qed

\end{document}